\documentclass[12pt]{amsart}
\usepackage[cp1251]{inputenc}
\usepackage[english,russian]{babel}
\usepackage{amsmath}
\usepackage{amssymb}
\usepackage{amsfonts}
\usepackage{graphicx}

\newtheorem{lemma}{Лемма}
\newtheorem{theorem}{Теорема}

\newtheorem{problem}{Задача}

\begin{document}

 УДК{517.984.50}

\title[Оператор градиент дивергенции]
{Оператор градиент дивергенции в пространствах Соболева II}
\author{\bf {Р.С.~Сакс}}


\maketitle
 Аннотация. {\small Автор изучает структуру пространства
  $\mathbf{L}_{2}(G)$ вектор-функций, квадратично интегрируемых
      по ограниченной   области $G$ трех-мерного    пространства
     с гладкой границей   $\Gamma$,
    и роль операторов   градиента дивергенции
     и ротора в построении базисов в его ортогональных подпро-странствах
     ${\mathcal{{A}}}$ и  ${\mathcal{{B}}}$.

В  ${\mathcal{{A}}}$ и  ${\mathcal{{B}}}$ выделяются
подпространства ${\mathcal{{A}}_{\gamma}}(G)\subset{\mathcal{{A}}}$
и $\mathbf{V}^{0}(G)\subset {\mathcal{{B}}}$.
   Оказывается, что операторы  $\nabla\mathbf{div}$ и $\mathrm{rot}$
продолжаются   в эти подпространства, их расширения $\mathcal{N}_d$
\,и $S$ являются самосопряженными и обратимыми, а их обратные
операторы $\mathcal{N}_{d}^{-1}$   и $S^{-1}$ -- вполне
непрерывными.

В каждом из этих подпространств  мы строим ортонормированный базис.
Объединяя эти базисы,  получаем полный ортонормированный базис
объемлющего пространства $\mathbf{L}_{2}(G)$, составленный из
собственных функций операторов градиента дивергенции и ротора.
 В  шаре  $B$,  базисные функции  определяются элементарными функциями.

Определены пространства $\mathcal{A}^{s}_{\mathcal{K}}(B)$.
Доказано, что условие $\mathbf{v}\in
\mathcal{A}^{s}_{\mathcal{K}}(B)$ необходимо и достаточно для
сходимости ее ряда Фурье (по собственным функциям градиента
дивергенции) в норме пространства Соболева $\mathbf{H}^{s}(B)$.

 Используя ряды Фурье
  функций $\mathbf{f}$ и $\mathbf{u}$,
автор исследует разрешимость (в пространствах $\mathbf{H}^{s}(G)$)
краевой задачи для системы
$\nabla\mathbf{div}\mathbf{u}+\lambda\mathbf{u}=\mathbf{f}$ в $G$,
 $\mathbf{n}\cdot\mathbf{u}|_{\Gamma}=g$ на границе,
при условии$\lambda\neq 0$.

В шаре $B$ краевая задача:
$\nabla\mathbf{div}\mathbf{u}+\lambda\mathbf{u}=\mathbf{f}$,
  $\mathbf{n}\cdot\mathbf{u}|_S=0$,
 решена полностью и для любых $\lambda$.
Доказано, что при  $\lambda\,\overline{\in}\, Sp\,
(-\nabla\mathrm{div})$,
 оператор $\nabla\text{div}+\lambda I$ осуществляет гомеоморфизм
 пространств  $\mathbf{{H}}^2_{\gamma\delta\gamma}(B)$
 и  $\mathbf{{F}}_{\gamma}^0(B)$.

Эта работа есть продолжение работы [17Ndi] "Оператор $\mathcal{N}_d$
в подпрост-ранствах $\mathbf{L}_{2}(G)$"(см.
arXiv:1710.06428v1[math.AP] 17 Oct2017)
 .{\footnote {\small Свойства оператора $S$ автор
изучал в работе [17S]: "Ряды Фурье оператора ротор и пространства
Соболева".}}

\section {Введение и основные результаты}

\subsection {Основные пространства} В статье мы рассматриваем линейные
пространства над полем $\mathbb{C}$ комплексных чисел. Через
$\mathbf{L}_{2}(G)$ обозначаем пространство Лебега вектор-функций,
квадратично интегрируемых в $G$    с внутренним произведением
$(\mathbf {u},\mathbf {v})= \int_G \,\mathbf
{u}\cdot\overline{\mathbf {v}}\,d\,\mathbf {x}$  и нормой
$\|\mathbf{u}\|~= (\mathbf {u},\mathbf {u})^{1/2}$.
  Пространство  Соболева, состоящее из вектор-функций,
    принадле-жащих $\mathbf{L}_{2}(G)$ вместе с обобщенными производными
  до порядка $s\geq 0$,   обозначается через
$\mathbf{H}^{s}(G)$, $\|\mathbf {f}\|_s$ -норма его элемента
$\mathbf {f}$;   замыкание в $\mathbf{H}^{s}(G)$ пространства
$\mathcal{C}^{\infty}_0(G)$ обозначается через
$\mathbf{H}^{s}_0(G)$. Пространство Соболева отрицательного порядка
$\mathbf{H}^{-s}(G)$ двойственно к $\mathbf{H}^{s}_0(G)$.
 {\footnote {\small см.
 $W_p^{(l)}(\Omega)$  в $\S 3$ гл. 4 \cite{sob},
    $H^k(Q)$ в $\S 4$ гл. 3 \cite{mi} и   гл.1 в    \cite{mag}.}}

  В области $G$  с гладкой границей  $\Gamma$ в каждой точке $y\in\Gamma$
  определена нормаль   $\mathbf {n}(y)$  к $\Gamma$.\,
Вектор-функция     $\mathbf {u}$ из $\mathbf{H}^{s+1}(G)$ имеет след
$ \gamma(\mathbf {n}\cdot\mathbf {u})$ на $\Gamma$ ее нормальной
компоненты, который принадлежит пространству Соболева-Слободетцкого
$\mathbf{H}^{s+1/2}(G)$, $|\gamma(\mathbf {n}\cdot\mathbf
{u})|_{s+1/2}$- его норма.

Пусть функция $h\in {H}^{1}(G)$, а $\mathbf{u}=\nabla h$ - ее
градиент. Через ${\mathcal{{A}}}(G)$ обозначим подпространство
$\{\nabla h, h\in H^1(G)\}$ в $\mathbf{L}_{2}(G)$, а через
${\mathcal{{B}}}(G)$ -  его ортогональное дополнение. Так что
\begin{equation}\label{wor  1}\mathbf{L}_{2}(G)=
{\mathcal{{A}}}(G)\oplus{\mathcal{{B}}}(G).\end{equation} {\small
Это разложение я взял из статьи Z.Yoshida и Y.Giga
  \cite{giyo}. Они ссылаются на работу Г.Вейля  \cite{hw}, который
          в   Теорема II
  утверждает, что $\mathfrak{F}_0=\mathfrak{G}+\mathfrak{F}'$, где
  $\mathfrak{F}_0={L}_2(G)$,   $\mathfrak{G}$ есть
  замыкание в норме ${L}_2$ градиентов функций $\psi\in
\Gamma$, $\Gamma=\mathcal{C}_0^1(G)$, а $\mathfrak{F}'$-множество
соленоидальных элементов в $\mathfrak{F}_0$.}

{\small Разложение $\mathbf{L}_{2}(\Omega)$ вводили и использовали
С.Л.Соболев \cite{so54},\,О.А.Лады-женская \cite{lad}, а также Кочин
Н.Е., Кибель И.А., Розе Н.В.\cite{kkr},\,\,Э.Быховский и Н.Смирнов
\cite{bs} и  К.Фридрихс \cite{fri}.}

Мы будем придерживаться разложения \eqref{wor  1}.

 Если граница
области $\Omega$ имеет положительный род $\rho$, то $ \mathcal{B}$
содержит в себе конечномерное подпространство
\begin{equation}\label{bh  1} \mathcal{B}_H=\{\mathbf{u}\in\mathbf{L}_{2}(\Omega):
\mathrm{rot}\,\mathbf{u}=0,\,\,\mathrm{div}\,\mathbf{u}=0,
\,\,\gamma(\mathbf{n}\cdot \mathbf{u})=0 \}. \end{equation} Его
размерность
 равна $\rho$ \cite{boso}, а базисные функции
 $h_j\in \mathcal{C}^\infty(\Omega)$\,\cite{hw}.
Гладкость обобщенных рещений системы
  \eqref{bh  1}\,
  доказал  Г.Вейль \cite{hw}.

   Ортогональное дополнение $\mathcal{B}_H$ в
$\mathcal{B}(\Omega)$,  обозначим через $\mathbf {V} ^{0} (\Omega)$.
Значит,
\begin{equation}\label{bhr  1}{\mathcal{B}(\Omega)=
\mathbf {V} ^{0} (\Omega)\oplus\mathcal{B}_H} (\Omega).
\end{equation}
{ \small {В \cite{giyo}:\,
${L}_{\sigma}^2(\Omega)={L}_{\Sigma}^2(\Omega)\oplus{L}_{H}^2(\Omega)$.
Символ $L$ перегружен. Мы изменили авторские обозначения пространств
${L}_{\Sigma}^2(\Omega)$
 и ${L}_{H}^2(\Omega)$  на
   $\mathbf {V}^{0} (\Omega)$ и $\mathcal{B}_{H} (\Omega)$. }}

   Индексом $\gamma$  будем снабжать пространства вектор-функций
        $\mathbf {u}$, нормальные компоненты которых
             имеют на границе $\Gamma$ нулевой    след
     $ \gamma(\mathbf {n}\cdot\mathbf {u})$:
     \begin{equation}\label{ag  1}\mathcal{A}_{\gamma} =
     \{\nabla h, h\in H^1(\Omega),  \gamma(\mathbf {n}\cdot {\nabla\,h})=0
     \}.\end{equation}
\subsection {Свойства операторов ротор и  градиент дивергенции}
Операторы градиент, ротор и   дивергенция определяются в трехмерном
векторном анализе \cite{zo}. Им соответствует оператор $d$ внешнего
диф-ференцирования на формах $\omega^k$ степени $k=0,1$ и 2.
Соотношения $dd\omega^k=0$ при $k=0,1$ имеют вид
$\mathrm{rot}\,\nabla h=0$ и
$\mathrm{div}\,\mathrm{rot}\,\mathbf{u}=0$.

Формулы $\mathbf{u}\cdot\nabla
h+h\mathrm{div}\mathbf{u}=\mathrm{div}(h \,\mathbf{u}) $,\,
$\mathbf{u}\cdot\mathrm{rot}\,\mathbf{v}-
\mathrm{rot}\,\mathbf{u}\cdot\mathbf{v}=\mathrm{div}[\mathbf{v},\mathbf{u}]$,
где $[\mathbf{v},\mathbf{u}]$ - векторное произведение, и
интегрирование по области $G$ используются при определении
операторов $\mathrm{div}\,\mathbf{u}$  и $\mathbf{rot}\mathbf{u}$ в
$\mathbf{L}_{2}(G)$   и $\mathcal{D}'(G)$.

 Оператор Лапласа
выражается через $\mathrm{rot}\,\mathrm{rot}$ и
$\nabla\,\mathrm{div}$:
\begin{equation}\label{dd  1}-\Delta \mathbf{v} =\mathrm{rot}\, \mathrm{rot}\,\mathbf{v}-
  \nabla \mathrm{div} \mathbf{v}.\end{equation}
  Оператор Лапласа эллиптичен, а операторы $\mathrm{rot}$ и
  $\nabla\,\mathrm{div}$ не являются таковыми \cite{so71}. Они
вырождены в том смысле, что   $\mathrm rot\,\mathbf {u}=0$ при
$\mathbf{u}\in~ \mathcal{A}(G)$, $\nabla\mathrm div\,\mathbf {v}=0$
при $\mathbf{v}\in \mathcal{B}(G)$   в 
$\mathcal{D}'(G)$, и являются аннуляторами друг друга.

При $\lambda\neq 0$ операторы $\mathrm{rot} +\lambda\, I$ на
$\mathcal{A}(G)$ и $\nabla\mathrm{div} +\lambda\, I$ на
$\mathcal{B}(G)$ сводятся к
  умножению на $\lambda$, а на ортогональных подпространствах,
   при 
   условии
   $ \gamma(\mathbf {n}\cdot\mathbf {u})=0$,--
    разрешимы по Фредгольму.

 Они принадлежат классу Б.Вайберга и В.Грушина \cite{vagr}
операторов, приводимых  к эллиптическим:  их расширения являются
эллиптическими (формально переопределенными) операторами. Краевые
задачи для этих расширений с условием
 $ \gamma(\mathbf {n}\cdot\mathbf {u})=g$ - эллиптичны по теории
 Солонникова  \cite{so71}. В результате мы получаем априорные оценки
 их решений (см. п.1.5)

  С другой стороны,  операторы  градиент дивергенции и ротор
допускают самосопряженные расширения в пространства
 ${\mathcal{{A}}_{\gamma}}(G)$ и $\mathbf{V}^{0}(G)$.

А именно, оператор
 $\mathcal{N}_d:\mathcal{A}_{\gamma}(G)
 \longrightarrow \mathcal{A}_{\gamma}(G)$ , совпадающий с
 $\nabla \mathrm{div}\mathbf {u}$
 на  $\mathcal{A}_ {\gamma}^{2}=\{\mathbf {u}\in\mathcal{A} _ {\gamma}(G):
 \nabla \mathrm{div}\mathbf {u}\in\mathcal{A} _ {\gamma}(G)\}$,  самосопряжен и
 его обратный оператор
 $\mathcal{N}^{-1}_d:\mathcal{A} _ {\gamma}\rightarrow \mathcal{A}_ {\gamma}^{2}$
  вполне непрерывен ( п.2.3).

 Соответственно, оператор \quad $S:\mathbf{V}^{0}(G)\longrightarrow
\mathbf{V}^{0}(G)$ с областью определения \quad
$\mathbf{W}^{1}=\{\mathbf {u}\in\mathbf {V}^0(G):\mathrm{rot}\mathbf
{u}\in\mathbf {V}^0(G)\}$,  совпадающий с $\mathrm{rot}\,\mathbf{u}$
при  $\mathbf {u}\in\mathbf{W}^{1}$,
 является самосопряженным и его обратный
 $S^{-1}$- вполне непрерывен \cite{giyo}.


Пусть $\mathbf{u}_{j}^{\pm}(\mathbf{x})$ --собственные функции
ротора, отвечающие ненулевым собственным значениям
$\pm{\lambda}_j\in {\Lambda}\subset\mathbb{R}$:
 $\mathrm{curl}\,\mathbf{u}_{j}^{\pm}={\pm}\lambda_j\mathbf{u}_{j}^{\pm}$,
$\|\mathbf{u}_{j}^{\pm}\|=1$, и они образуют ортонормированный базис
в $\mathcal{B}(G)$. Проекции векторов $\mathbf{f}$ и $S\mathbf{f}$
из $\mathbf{L}_{2}(G)$ на $\mathcal{B}(G)$ имеют вид
 \begin{equation}\label{fb__1_}
\mathbf{f}_{\mathcal{B}}=\sum_{\pm\lambda_j\in
\Lambda}[(\mathbf{f},\mathbf{u}_{j}^{+})
\mathbf{u}_{j}^{+}(\mathbf{x})+(\mathbf{f},\mathbf{u}_{j}^{-})
\mathbf{u}_{j}^{-}(\mathbf{x})], \quad
\mathbf{f}(\mathbf{x})\in\mathbf {W} ^ {1}(G),
\end{equation}
\begin{equation}\label{sfb__2_}
S\mathbf{f}_{\mathcal{B}}=\sum_{\pm\lambda_j\in
\Lambda}\,\lambda_j[(\mathbf{f},\mathbf{u}_{j}^{+})
\mathbf{u}_{j}^{+}(\mathbf{x})-(\mathbf{f},\mathbf{u}_{j}^{-})
\mathbf{u}_{j}^{-}(\mathbf{x})].
\end{equation}
 Оператор $\mathcal{N}_d$  имеет в $\mathcal{A}_{\gamma}$
полную систему собственных функций $\mathbf{q}_{j}(\mathbf{x})$,
отвечающих ненулевым собственным значениям $\mu_j\in M\subset
\mathbb{R}$:
\begin{equation}\label{fa  1}-\nabla\mathrm{div}\,\mathbf{q}_{j}(\mathbf{x})=
 \mu_j\,\mathbf{q}_{j}(\mathbf{x}), \quad \|\mathbf{q}_{j}\|=1,
  \quad \text{и}\quad
\mathbf{f}_\mathcal{A}=\sum_{\mu_j\in M}
(\mathbf{f},\mathbf{q}_{j})\mathbf{q}_{j} (\mathbf{x}),
\end{equation}
\begin{equation}\label{gdf  2}\mathcal{N}_d\,\mathbf{f}_\mathcal{A}=-
\sum_{\mu_j\in M}\,\mu_j (\mathbf{f},\mathbf{q}_{j})\mathbf{q}_{j}
(\mathbf{x}), \quad \mathbf{f}(x)\in
\mathcal{A}^2_{\gamma}(G).\end{equation}
 {\small В  шаре $B$ радиуса $R$ собственные функции
 $\mathbf{u}^{\pm}_{\kappa}$ ротора, отвечающие
  ненулевым собственным значениям
  $\pm\lambda_{\kappa}=\pm\rho_{n,m}/R$ и собственные функции
  $\mathbf{q}_{\kappa}$
градиента дивергенции с собственными значениями $\nu_{\kappa}^2$,
$\nu_{\kappa}=\alpha_{n,m}/R,$ выражаются  явными формулами
 \cite{saUMJ13} и 
    \[\mathbf{rot}\,\mathbf{u}^{\pm}_{\kappa}=\pm\lambda_{\kappa}\,
\mathbf{u}^{\pm}_{\kappa}, \quad
\gamma\mathbf{n}\cdot\mathbf{u}^{\pm}_{\kappa}=0;\quad
\mathbf{rot}\,\mathbf{q}_{\kappa}=0
 \quad \kappa=(n,m,k),\,\]
\[ \nabla\,div\,\mathbf{u}^{\pm}_{\kappa}=0, \quad
-\nabla\,div\,\mathbf{q}_{\kappa}=\nu_{\kappa}^2\mathbf{q}_{\kappa},
\quad \gamma\mathbf{n}\cdot\mathbf{q}_{\kappa}=0, \quad
 \quad\,\, |k|\leq n ,\] где
числа $\pm\rho_{n,m}$ и $\alpha_{n,m}$ - нули функций $\psi_n$ и их
производных $\psi_n'$, а
 \begin{equation}
\label{psi__1_}\psi_n(z)=(-z)^n\left(\frac{d}{zdz}\right)^n\frac{\sin
z}z, \quad n\geq 0,\,\,m\in \mathbb{N}.\end{equation}

Cобственные функции каждого из операторов взаимно ортогональны и их
совокупная система полна в ${\mathbf{{L}}_{2}}(B)$ \cite{saUMJ13}.}

 \subsection{ Структура работы и основные результаты}

В $\S\, 2$  в ограниченной  области $G$ с гладкой границей $\Gamma$
изучается разрешимость краевой задачи
\begin{equation}
\label{ndi__1_}\nabla\text{div}\mathbf{u}+\lambda\,
\mathbf{u}=\mathbf{f}(\mathbf{x}), \quad
 \mathbf{x}\in G, \quad \mathbf{n}\cdot \mathbf{u}{{|}_{\Gamma }}=g
 \end{equation}
в пространствах Соболева $\mathbf{H}^{s+2}(G)$,  где число $s\geq 0$
целое,  $\mathbf{f}$ и $g$ заданы. Задаче ставится в соответствие
ограниченный  оператор $\mathbb{B}$.

Утверждается, что эта задача обобщенно эллиптична при $\lambda\neq
0$, то-есть она приводится к эллиптической задаче по Солонникову
\cite{so71}. Откуда вытекает Теорема 1:
 конечномерность ядра оператора $\mathbb{B}$ задачи в
пространстве Соболева $\mathbf{H}^{s+2}(G)$ и априорная оценка:
\begin{equation}
\label{ond__1_} C_s\|\mathbf{u}\|_{s+2}
\leq|\lambda|\|\mathrm{rot}^2\,\mathbf{u}\|_{s}+
\|\nabla\mathrm{div}\,\mathbf{u}\|_{s}+
|\gamma({\mathbf{n}}\cdot\mathbf{u})|_{s+3/2}+ \|\mathbf{u}\|_{s}.
\end{equation}
В п.2.3 мы изучаем оператор $\nabla\mathrm{div} +\lambda \mathbf{I}
$ в ортогональных подпространствах $\mathcal{A}$ и $\mathcal{B}$ в
$\mathbf{L}_2(G)$. На  $\mathcal{B}$
  оператор $\nabla\mathrm{div} \mathbf{u}
+\lambda \mathbf{u} $ сводится к
 $\lambda\mathbf{u}$. На
 $\mathcal{A}_{\gamma}$  он продолжается
 как самосопряженный оператор $\mathcal{N}_d + \lambda\,I$.
 Найдены необходимые и достаточные условия его
обратимости (Теорема 2b).

В  $\S\, 3$ установлены соотношения между решениями  спектральной
задачи~1 для градиента дивергенции (в
  области с гладкой границей)   и спектральной
задачи Неймана для скалярного оператора Лапласа (Лемма 1).

 В шаре ее решения
вычислены явно\,\cite{vla}. В результате мы получаем формулы
\eqref{qom 1} собственных функций $\mathbf{q}_{\kappa}(\mathbf{x})$
 градиента дивергенции.

Определены пространства $\mathcal{A}^{s}_{\mathcal{K}}(B)$.
Доказано, что условие $\mathbf{v}\in
\mathcal{A}^{s}_{\mathcal{K}}(B)$ необходимо и достаточно для
сходимости ее ряда Фурье (по собственным функциям градиента
дивергенции) в норме пространства Соболева $\mathbf{H}^{s}(B)$
порядка $s>0$ (Теорема 4).

В   $\S 4$ мы показываем, что
 собственные вектор-функции  градиента дивер-генции  и ротора
$$\{\mathbf{q}_{i}(\mathbf{x}), \,\mathbf{h}_{j}(\mathbf{x}),\,
 \mathbf{q}_{k}^{+}(\mathbf{x}),\, \mathbf{q}_{k}^{-}
 (\mathbf{x})\}\quad
\mu_i\in M,\,\, j\in [1,\rho],\, \,\pm\lambda_k\in \Lambda,$$
взаимно ортогональны и образует в  ${\mathbf{{L}}_{2}}(G)$
ортонормированный базис. Приводится   разложение
 векторного поля $\mathbf{f}$ из
${\mathbf{L}_{2}}(B)$ на потенциальное и соленоидальное   поле:
 $\mathbf{f}(\mathbf{x})=\mathbf{f}_{\mathcal{A}}(\mathbf{x})+
 \mathbf{f}_{\mathcal{B}}(\mathbf{x})$ 
 и изучаются свойства рядов $\mathbf{f}_{\mathcal{A}}$ и
 $\mathbf{f}_{\mathcal{B}}$
в пространствах  $\mathbf{{F}}_{\gamma}^0(B)$
 и  $\mathbf{{E}}_{\gamma}^0(B)$.

В  $\S\, 5$   методом Фурье решена краевая задача
$\nabla\mathbf{div}\mathbf{v}+\lambda\mathbf{v}=\mathbf{f}$ в шаре
$B$ с условием $\gamma\mathbf{n}\cdot\mathbf{v}=0$   при любых
$\lambda$ (Теорема 6).

Оператор $\nabla\mathbf{div}+\lambda \,I$ отображает пространство
  $\mathbf{{H}}^2_{\gamma\delta\gamma}(B)$ в
  $\mathbf{{F}}_{\gamma}^0(B)$ непрерывно.
  Доказано, что при  $\lambda\,\overline{\in}\, Sp\,
(-\nabla\mathrm{div})$
 оператор $\nabla\text{div}+\lambda I$ осуществляет гомеоморфизм
 пространств  $\mathbf{{H}}^2_{\gamma\delta\gamma}(B)$
 и  $\mathbf{{F}}_{\gamma}^0(B)$ \, (Лемма 2).

\section{ Градиент дивергенции
 в ограниченной области}
\subsection{Краевая задача:}
в ограниченной  области $G$ с гладкой границей $\Gamma$ заданы
векторная и скалярная функции $\mathbf{f}$ и ${g}$, найти
вектор-функцию $\mathbf{u}$, такую что
\begin{equation}
\label{grd__2_}
 \nabla\, \text{div}\mathbf{u}+\lambda\, \mathbf{u}=\mathbf{f}
 \quad  \text{в}\quad G,
\quad \mathbf{n}\cdot\mathbf{u}|_\Gamma=g,
\end{equation}
  Эта задача при  $\lambda\ne 0$  является обобщенно
   эллиптической.

 Это означает, что  при  $\lambda\ne 0$ система \eqref{grd__2_}
  принадлежит классу (RNS1) Вайнберга и Грушина \cite{vagr},
  а краевая задача для расширенной системы
  эллиптична по Солонникову \cite{so71}.

   Действительно,
  оператор $\nabla\, \text{div}+\lambda\text{I}$ второго порядка таков, что \newline
   а)  {\it ранг его символической матрицы
$\nabla\, \text{div}(i\xi)$ постоянный (и равен единице) при любых
$\xi\neq0$},
\newline
 б){\it оператор $\nabla\,\text{div}$  имеет  левый
 аннулятор:} это $\text{rot}$, так как $ \text{rot}\,\nabla\, \text{div} \mathbf{u}=0$
 для любой $\mathbf{u}\in \mathcal{D}'(G)$,\newline
 в){\it расширенная система:
  \begin{equation}
\label{rodi__2_} \nabla\, \text{div}\,\mathbf{u}+\lambda
\mathbf{u}=\mathbf{f},\quad \lambda \text{rot}\,\mathbf{
u}=\text{rot}\,\mathbf{f}\end{equation} является эллиптической (по
Даглису-Ниренбергу)}. Ранг ее символической матрицы
$\nabla\,\text{div}(i\xi )$// $\lambda\text{rot}(i\xi) $ равен  трем
при всех $\xi \neq 0$. {\footnote {\small При выборе порядков $s_k$
и $t_j$ для его строк  и столбцов:  $s_k=0$ для $k=1,2,3$ и $s_k=-1$
для $k=4,5,6$; и при $t_j=2$ для $j=1,2,3$\, \cite{sa75}.}}

Эта система переопределена, если вместо $\text{rot}\,\mathbf{f}$
стоит 
    любая $\mathbf{F}$.

 Оказывается, что 
 краевая задача  $\gamma\mathbf{n}\cdot\mathbf{u}=g$ системы
 \eqref{rodi__2_},
  эллиптична по Солонникову \cite{so71}. Это означает что

1)  система \eqref{rodi__2_} эллиптична,

2) граничный оператор $\gamma\mathbf{n}\cdot\mathbf{u}$ "накрывает"
\,оператор системы \eqref{rodi__2_}.

Первое условие выполнено.

 Пусть  $\tau$ и $\mathbf{n}$- касательный
и нормальный векторы к $\Gamma$ в точке $y\in \Gamma$ и
$|\mathbf{n}|=1$.\quad
  Второе условие выполняется, если

  $2^0)$   однородная система линейных
дифференциальных уравнений ( на полуоси\, $z\geq 0$ с параметром
$|\tau| > 0$):
\begin{equation}\label{cdz  3}
\lambda\text{rot}(i\tau+\mathbf{n} d/dz ) \mathbf{v}=0,\quad
(\nabla\text{div})(i\tau+\mathbf{n} d/dz )\mathbf{v}=0 \,\, 
\end{equation} при условиях:
$\mathbf{n}\cdot \mathbf{v}|_{ z=0}=0$ и $\mathbf{v}\rightarrow 0$
при $z\rightarrow + \infty$ имеет только тривиальное решение
$\mathbf{v}(y,\tau; z )$.

   Это доказано в 
   работе   автора в  [17Ndi] "Оператор  градиент
дивергенции  в подпространствах $\mathbf{L}_{2}(G)$" (см. ArXiv.org
:1710.06428v1[math.AP] 17 Oct2017).
 Здесь же мы изложим другие основные результаты.

\subsection{Оператор задачи \eqref{grd__2_} в пространствах Соболева}
Пусть $\mathbf{u}$ принадлежит пространству
${\mathbf{H}^{s+2}}(G)$,\,  $s\geq  0$.
 Тогда
$\nabla\text{div} \mathbf{u}$ и $\text{rot}^{2}\mathbf{u}$
принадлежат $\mathbf{H}^{s}(G)$.
 Поэтому вектор-функция $\mathbf{f}:=\nabla\text{div} \mathbf{u}+ \lambda
\mathbf{u}$ принадлежит пространству
\begin{equation}\label{pr  2} {\bf{F}^{s}}(G)=\{\mathbf{f}\in
{\mathbf{H}^{s}}: \text{rot}^2\,\mathbf{f}\in
 \mathbf{H}^{s}\},\quad  
 \|\mathbf{v}\|_{\mathbf{F}^{s}}=(
 \|\mathbf{v}\|^2_{s}+
 \|\text{rot}^2\mathbf{v}\|^2_{s})^{1/2}.\end{equation} Функция
$g:=\gamma({\mathbf{n}}\cdot\mathbf{u})$
принадлежит пространству
$H^{s+3/2}(\Gamma)$. Следовательно, при $\lambda\neq 0$ задаче
соответствует ограниченный оператор
\begin{equation}\label{op  2} \mathbb{B}\mathbf{u}\equiv \begin{matrix}
\nabla\text{div}\,\mathbf{u}+\lambda\,\mathbf{u} \\
\gamma({\mathbf{n}}\cdot\mathbf{u})\end{matrix}:
\bf{H}^{s+2}(G)\rightarrow
\begin{matrix}\bf{F}^{s}(G)\\ H^{s+3/2}(\Gamma)\end{matrix}.
\end{equation}
Согласно Теореме 1.1 из работы Солонникова \cite{so71}, о
переопределенных эллиптических краевых задачах  в ограниченной
области $G$   с гладкой границей $\Gamma\in \mathcal{C}^{s+2} $,
обобщенно эллиптический оператор \eqref{op 2} имеет левый
регуляризатор:
 ограниченный оператор  $\mathbb{B}^L$ такой, что
$\mathbb{B}^L\mathbb{B}=\mathbb{I}+\mathbb{T}$, где $\mathbb{I}$ -
единичный, а $\mathbb{T}$ - вполне непрерывный операторы, и
существует постоянная $C_s >0$ такая, что выполняется априорная
оценка:
\begin{equation} \label{apgro s} C_s\|\mathbf{u}\|_{s+2}
\leq|\lambda|\,\|\mathrm{rot}^2\,\mathbf{u}\|_{s}+
\|\nabla\mathrm{div}\,\mathbf{u}\|_{s}+
|\gamma({\mathbf{n}}\cdot\mathbf{u})|_{s+3/2}+ \|\mathbf{u}\|_{s}.\end{equation}
 { \small И такая же
 оценка в  пространствах Соболева    $W^{s+2}_p(G)$, $
p>1$.}
 Значит, имеет место
\begin{theorem}
Оператор $\mathbb{B}$ в пространствах (\ref{op  2}) имеет левый
регуля-ризатор. Его ядро $\mathcal{M}$ конечномерно и выполняется
 оценка \eqref{apgro s}.
 \end{theorem}

Из этой теоремы и  оценки следует, что при $\lambda \ne 0$

 a){\it число линейно независимых решений однородной
задачи \eqref{grd__2_} конечно,}

 b){\it любое ее обобщенное решение 
 бесконечно дифференцируемо вплоть до
границы, если граница области бесконечно дифференцируема.}

\subsection{Оператор $\nabla\mathrm{div}+\lambda I$ в подпространствах }
На подпространстве $\mathcal{B}$ в $\mathbf{L}_2(G)$, ортогональном
подпространству
  $\mathcal{A}$,
 оператор $\nabla\, \text{div}\mathbf{u} +\lambda \mathbf{u} $
является  оператором умножения: $\lambda \mathbf{u} $.

 Пространство
$\mathcal{A}_{\gamma}=\{\mathbf{u}=\nabla\, h:\,\, h\in H^1(G),\quad
(\mathbf{n}\cdot\mathbf{u})|_{\Gamma}=0 \}$ плотно в $\mathcal{A}$,
так как функции из $\mathcal{C}_0^\infty \cap \mathcal{A}_{\gamma}$
плотны в $\mathbf{L}_2(G)$. Пространство
  \begin{equation}
\label{sa__1_}\mathcal{A}^2_{\gamma}(G)=\{ \mathbf{v}\in
\mathcal{A}_{\gamma}: \nabla\mathrm{div} \mathbf{v}\in
\mathcal{A}_{\gamma}\}.\end{equation}  плотно в
$\mathcal{A}_{\gamma}$ и
 содержится в $\mathbf{H}^2(G)$ в силу оценки \eqref{apgro s}.

Введем оператор $\mathcal{N}_d:\mathcal{A}_{\gamma}
 \rightarrow \mathcal{A}_{\gamma}$ с  областью
определения   $\mathcal{A}^2_{\gamma}(G)$, 
который совпадает с $\nabla\mathrm{div}\,\mathbf{v}$ при
$\mathbf{v}\in \mathcal{A}^2_{\gamma}(G)$.

Оператор $\mathcal{N}_d+\lambda I:\mathcal{A}_{\gamma}\rightarrow
\mathcal{A}_{\gamma}$ является самосопряженным (эрмитовым
\,\cite{vla}).
 Действительно, согласно формуле Гаусса-Остроградсного
\begin{equation}
\label{gri__2_} \int_G(\nabla\mathrm{div}\mathbf{u}+\lambda
\mathbf{u})\cdot\,\mathbf{v} dx=
\int_G\mathbf{u}\cdot\,(\nabla\mathrm{div}\mathbf{v}+\lambda
\mathbf{v})dx+\end{equation}
\[\int_{\Gamma}[(\mathbf{n}\cdot\mathbf{v})\mathrm{div}\mathbf{u}+
(\mathbf{n}\cdot\mathbf{u})\mathrm{div}\mathbf{v}]|_\Gamma \, d S.\]

Если вектор-функции $\mathbf{u}$ и $\mathbf{v}$ принадлежат
$\mathcal{A}^2_{\gamma}(G)$, то граничные интегралы пропадают,
остальные интегралы сходятся. Следовательно,
\begin{equation}
\label{gri__2_} ((\nabla\mathrm{div}\mathbf{u}+\lambda
\mathbf{u}),\mathbf{v})=
(\mathbf{u},(\nabla\mathrm{div}\mathbf{v}+\lambda \mathbf{v})) \quad
\text{в}\quad \mathbf{L}_2(G).\end{equation}

{\it Область определения $\mathcal{A}^2_{\gamma}(G)$ оператора
$\mathcal{N}_d$ содержится в $\mathbf{H}^2(G)$ и плотна в
$\mathcal{A}_{\gamma}$, а  область его значений совпадает с
$\mathcal{A}_{\gamma}$} [17Ndi].

Пространство $\mathcal{A}_\gamma$ ортогонально ядру оператора
$\mathcal{N}_d$, он имеет единственный обратный
$\mathcal{N}_d^{-1}$, определенный на $\mathcal{A}_\gamma$. Оператор
$\mathcal{N}_d^{-1}:\mathcal{A}_{\gamma}\rightarrow
 \mathcal{A}^2_{\gamma}(G)$ имеет точечный спектр, который
 не содержит точек накопления кроме нуля
Следовательно,  спектр самосопряженного  оператора $\mathcal{N}_d$
точечный  и действительный, а система его собственных вектор-функций
ортогональна и полна в пространстве $\mathcal{A}_{\gamma}$. Каждому
собственному значению соответствует конечное число  собственных
вектор-функций.

Пусть $\mathbf{f}\in \mathcal{A}_{\gamma}(G)$, так как
 $(\mathcal{N}_d+\lambda\,I)\mathbf{f}\in \mathcal{A}_{\gamma}(G)$
то
\begin{equation}
\label{sp__2_} (\mathcal{N}_d+\lambda\,I)\mathbf{f}=\sum_{\mu_j\in
{M}}[(\lambda-\mu_j)(\mathbf{f},\mathbf{q}_{j})\mathbf{q}_{j} ]
\end{equation}и ряд сходится в $\mathbf{L}_{2}(G)$.
Если $\lambda=\mu_{j_0}$,
 то соответствующее слагаемое в этом ряду исчезает.

Если элемент $(\mathcal{N}_d+\lambda\,I)^{-1}\mathbf{f}\in
\mathcal{A}_{\gamma}(G)$, то
\begin{equation}
\label{sp__4_}
(\mathcal{N}_d+\lambda\,I)^{-1}\mathbf{f}=\sum_{\mu_j\in
M}[(\lambda-\mu_j)^{-1}(\mathbf{f},\mathbf{q}_{j})\mathbf{q}_{j} ]
\end{equation} и ни одно из слагаемых этого ряда не обращается в
бесконечность. Это означает, что $(\mathbf{f},\mathbf{q}_{j})=0$ при
$\lambda=\mu_j=\mu_{j_0}$, то-есть функция $\mathbf{f}$ ортогональна
всем собственным функциям $\mathbf{q}_{j}(\mathbf{x})$ градиента
дивергенции, отвечающим собственному значению $\mu_{j_0}$.
 Итак, имеет место
\begin{theorem}
a). Оператор $\mathcal{N}_d:\mathcal{A}_{\gamma}\rightarrow
\mathcal{A}_{\gamma}$   является
 самосопряженным.
  Его спектр  $\sigma(\mathcal{N}_d )$
 точечный  и действительный.
  Семейство собственных функций $\mathbf{q}_{j}(x)$ оператора
  $\mathcal{N}_d$ образует полный ортонормированный базис в  пространстве
   $\mathcal{A}_{\gamma}$;
    разложение   $\mathbf{a}(x)\in{\mathcal {{A}} _{\gamma}} (G)$ имеет вид
\begin{equation}
\label{spr__4_}\mathbf{a}(x)=\sum_{\mu_j\in
M}(\mathbf{a},\mathbf{q}_{j}) \mathbf{q}_{j}(x),\quad
\|\mathbf{q}_{j}\|=1.\end{equation}
  b). Если $-\lambda$ не совпадает ни с одним из
   собственных   значений оператора $\mathcal{N}_d$, то
   оператор $\mathcal{N}_d+\lambda\,I:
    \mathcal{A}_{\gamma}\rightarrow\mathcal{A}_{\gamma}$ однозначно обратим,
    и его обратный   задается формулой   \eqref{sp__4_}.
Если $\lambda=\mu_{j_0}$, то он обратим тогда и только
  тогда,когда
  \begin{equation}
\label{urz _1_}\int_G \mathbf{f}\cdot \mathbf{q_j}\, dx=0\quad
\text{для}\,\,\forall \mathbf{q_j}: 
\mu_j=\mu_{j_0}.
\end{equation}
Ядро оператора $\mathcal{N}_d-\mu_{j_0}\,I$ определяется
собственными функциями $\mathbf{q_j}(\mathbf{x})$, собственные
значения которых равны $\mu_{j_0}$:
\begin{equation} \label{ker__1_}
Ker(\mathcal{N}_d-\mu_{j_0}\,I)= \sum_{\mu_j=\mu_{j_0}}
c_j\,\mathbf{q}_{j}(\mathbf{x}), \quad \text{для}\,\,\forall c_j \in
\mathcal{R}.\end{equation}
  \end{theorem}

\section {Связь между собственными функциями операторов $\nabla
div$ и Лапласа-Неймана}
\begin{problem} Найти  все ненулевые собственные значения $\mu
$ и  собственные вектор-функции $\mathbf{u}(\mathbf{x})$ в
${{\mathbf{L}}_{2}}(G)$ оператора градиент дивергенции такие, что
\begin{equation}  \label{gd   1}-\nabla \text{ div }\mathbf{u}=
\mu \mathbf{u}\quad  \text {в} \quad G,\quad \mathbf{n}\cdot
\mathbf{u}{{|}_{\Gamma }}=0,
\end{equation}
              где
$\mathbf{n}\cdot \mathbf{u}$ - проекция вектора $\mathbf{u}$ на
нормальный вектор $\mathbf{n}$.\end{problem}

Эта задача связана со спектральной задачей Неймана для скалярного
оператора Лапласа.

 \begin{problem} Найти все собственные значения
$\nu $ и собственные функции $g (\mathbf{x})$  оператора Лапласа
$-\Delta $ такие, что
      \begin{equation}  \label{gen   1}
              -\Delta g =\nu g \quad\text{в} \,\, G,\quad
 \mathbf{n}\cdot\nabla\,g|_{\Gamma }=0.
              \end{equation}\end{problem}

Эта задача является самосопряженной \cite{vla, mi}. Решения задач 3
и 4 принадлежат классу $\mathcal{C}^\infty(\overline{G})$, так как
$\Gamma\in \mathcal{C}^\infty$.

 Легко убедиться непосредственно, что
 \begin{lemma} Любому решению
$(\mu ,\mathbf{u})$ задачи 3 в области G соответствует
 решение $(\nu ,g )=(\mu, \text{div }\mathbf{u})$ задачи 4.
  Обратно, любому решению  $(\nu ,g )$ задачи 4
соответствует решение  $(\mu ,\mathbf{u})=(\nu, \nabla g)$ задачи
3.\end{lemma}
    \subsection {Явные решения спектральной задачи Лапласа-Неймана в шаре }
   Согласно книге \cite{vla} В.С.Владимирова

    {\it собственные значения
    оператора Лапласа-Неймана $\mathcal{N}_\Delta$ в шаре  $B$   равны
    $\nu _{n,m}^{2}$,  где   $\nu
     _{n,m}^{{}}={{\alpha }_{n,m}}/R$,    $n\ge 0$,   $m\in N$, а   числа
      ${{\alpha }_{n,m}}>0$ суть нули  функций
      ${{\psi }_{n}^{\prime }}(z)$,
     производных ${{\psi }_{n}}(z)$.
Соответствующие $\nu _{n,m}^{2}$ собственные функции $g _{\kappa
}^{{}}$ имеют вид:
 \begin{equation}  \label{lan   1}g _{\kappa }^{{}}(r,\theta ,\varphi )=c{{_{\kappa }^{{}}}^{{}}}
 {{\psi }_{n}}{{(\alpha _{n,m}^{{}}r/R)}^{{}}}Y_{n}^{k}(\theta ,\varphi
 ), \quad \kappa =(n, m, k),  \end{equation}
где   $c_{\kappa }^{{}}$-произвольные действительные
постоянные,\,$\kappa$- мультииндекс, \, $Y_{n}^{k}(\theta ,\varphi
)$ - действительные сферические функции,  $n\ge 0$,\, $|k|\le n, \,
m\in N$.}

Функции $g _{\kappa }^{{}}(x)$  принадлежат классу ${{C}^{\infty
}}(\overline{B})$  и при различных $\kappa$ ортогональны  в
${{L}_{2}}(B)$.   Система   функций $\{g _{\kappa }^{{}}\}$
полна в ${{L}_{2}}(B)$ \cite{mi}.
Нормируя их, получим  ортонормированный в ${{L}_{2}}(B)$ базис.


\subsection {Решение спектральной задачи 1 для $\nabla div$ в шаре}
Согласно лемме 1 вектор-функции ${{\mathbf{q}}_{\kappa }}(x)=
     \nabla {{g }_{\kappa }}(x)$
     являются решениями задачи 1 при ${\mu}_{n,m} ={\alpha }_{n,m}^2R^{-2}$ в
     ${{\mathbf{L}}_{2}}(B)$.
     Их компоненты $(q_r,q_\theta, q_\varphi)$ имеют вид
     \begin{equation}  \label{qom   1}\begin{array}{c}
     q _{r,\kappa }^{{}}(r,
     \theta ,\varphi )=c_{\kappa }(\alpha _{n,m}/R)
 {{\psi }_{n}^{{\prime}}}{{(\alpha _{n,m}^{{}}r/R)}^{{}}}Y_{n}^{k}
 (\theta ,\varphi ),\\
(q_{\varphi}+iq_{\theta})_{\kappa}=c_{\kappa }(1/r)
 {\psi }_{n}(\alpha _{n,m}r/R)\text{H}Y_{n}^{k}
 (\theta ,\varphi ),
  \end{array}\end{equation}
  \[\text{H}Y_{n}^{k}(\theta ,\varphi )=
  {{\left( {{\sin }^{-1}}{{\theta }^{_{{}}}}{{\partial }_{\varphi }}+
  i{{\partial }_{\theta }} \right)}^{{}}}Y_{n}^{k}(\theta ,\varphi ).\]
При $\kappa=(0,m,0)$ функция $Y_{0}^{0} (\theta ,\varphi )=1$,
$\text{H}Y_{0}^{0} =0$. Поэтому
\begin{equation}  \label{qomo   1}\begin{array}{l}
     q _{r,(0,m,0) }^{{}}(r)=c_{(0,m,0) }(\alpha _{0,m}/R)
 {{\psi }_{0}^{{\prime}}}{{(\alpha _{0,m}r/R)}},\\
(q_{\varphi}+iq_{\theta})_{(0,m,0)}=0.
  \end{array}\end{equation}

 Отметим, что $ {{\mathbf{q} }_{\kappa }}$ и $ {{\mathbf{q}
}_{{{\kappa }'}}}$ ортогональны при  ${\kappa }'\ne \kappa $.




\subsection{Сходимость ряда  по системе собственных функций
$\mathbf{q}_{\kappa}(\mathbf{x})$
 оператора $\nabla\mathrm{div}$ 
  в норме пространства Соболева $H^s(B)$} Определим подпространство
   $\mathbf{A}^s_\mathcal{K}(B)$ в $\mathcal{A}$
при $s\geq 1$:
\[ \mathbf{A}^s_\mathcal{K}(B)= \{\mathbf{f}\in
\mathcal{A}\cap\mathbf{H}^s(B): \mathbf{n}\cdot\mathbf{f}|_S=0,...,
 \mathbf{n}\cdot\ (\nabla\text{div})^{\sigma}\mathbf{f}|_S=0,\,\,
 \|\mathbf{f}\|_{\mathbf{A}_\mathcal{K}^s}=
\|\mathbf{f}\|_{\mathbf{H}^s}\},\]
 где $\sigma=[(s)/2]$.
  Имеет место
\begin{theorem}  Для того, чтобы $\mathbf{f}\in
\mathcal{A}$
 разлагалась в ряд Фурье
 \begin{equation} \label{arof 1}
\mathbf{f}(\mathbf{x})=\sum_{\kappa}
(\mathbf{f},\mathbf{q}_{\kappa})\mathbf{q}_{\kappa}(\mathbf{x})
\end{equation}
 по системе собственных вектор-функций $\mathbf{q}_{\kappa}(\mathbf{x})$
 оператора градиента дивергенции в шаре,
 сходящийся в норме
 пространства Соболева $\mathbf{H}^s(B)$, необходимо и достаточно,
  чтобы $\mathbf{f}$ принадлежала $\mathbf{A}^s_\mathcal{K}(B)$.

 Если $\mathbf{f}\in \mathbf{A}^s_\mathcal{K}(B)$,
то сходится ряд
\begin{equation} \label{arof 2}
\sum_{\kappa}{\nu}_{\kappa}^{2s}\,
|(\mathbf{f},\mathbf{q}_{\kappa})|^2 ,\quad
{\nu}_{\kappa}=({\alpha}_{n,m})/R
\end{equation} и существует такая положительная постоянная $C>0$, не
зависящая от $\mathbf{f}$, что
\begin{equation} \label{aorf 3}
\sum_{\kappa} {\nu}_{\kappa}^{2s}\,
|(\mathbf{f},\mathbf{q}_{\kappa})|^2 \leq
C\|\mathbf{f}\|^2_{\mathbf{H}^s(B)}.
\end{equation}
  Если $s\geq 2$, то любая вектор-функция  $\mathbf{f}$ из
  $\mathbf{A}^s_\mathcal{K}(B)$
разлагается в в ряд Фурье, сходящийся в пространстве $\mathbf{C}^{s-2}(\overline{B})$.%
\end{theorem}
Доказательство этой теоремы см. в [17Ndi].

{\bf Следствие.} {\it  Вектор-функция $f$ из
$\mathcal{A}\cap\mathbf{C}^{\infty}_0({B})$ разлагается в
 ряд Фурье \eqref{arof 1},
сходящийся в любом из пространств $\mathbf{C}^{k}(\overline{B})$,
$k\in \mathbb{N}$.}

\subsection{Скалярное произведение функций $\mathbf{f}$ и
$\mathbf{g}$ из $\mathcal{A}_{\gamma}$ в базисе из собственных
функций градиента дивергенции} Оно имеет вид:
\begin{equation} \label{spa 2}( \mathbf{f}, \mathbf{g})=
\sum_{\kappa, n\geq 0} \,
(\mathbf{f},\mathbf{q}_{\kappa})(\mathbf{g},\mathbf{q}_{\kappa})
\end{equation}
Если $\mathbf{f}$ и $\mathbf{g}$ принадлежат
$\mathbf{A}^2_\mathcal{K}(B)$, то равенства
\begin{equation} \label{sgd }(\nabla\bf{div}\, \mathbf{f}, \mathbf{g})=
(\mathbf{f},\nabla\text{div}\mathbf{g})= \sum_{\kappa, n\geq 0}
{\nu}_{\kappa}^2[(\mathbf{f},\mathbf{q}_{\kappa})
(\mathbf{g},\mathbf{q}_{\kappa}) ]
\end{equation}показывают, что оператор $\nabla\mathbf{div}$ является
самосопряженным в $\mathcal{A}_{\gamma}$.

\subsection { Сходимость ряда по собственным функциям
 ротора в норме пространства Соболева $\mathbf{H}^s(B)$, $s\geq 1$.}
 В [17S] доказана
  \begin{theorem}
 Для того, чтобы $\mathbf{f}\in \mathbf{V}^0(B)$
 разлагалась в ряд Фурье
 \begin{equation} \label{rof 1}
\mathbf{f}(\mathbf{x})=\sum_{\kappa, n>0}
((\mathbf{f},\mathbf{q}_{\kappa}^+)\mathbf{q}_{\kappa}^+(\mathbf{x})
+(\mathbf{f},\mathbf{q}_{\kappa}^-)\mathbf{q}_{\kappa}^-(\mathbf{x})),
\quad \|\mathbf{q}_{\kappa}^{\pm}\| =1,
\end{equation}
 по  собственным вектор-функциям $\mathbf{q}_{\kappa}^{\pm}(\mathbf{x})$
 ротора в шаре,
 сходящийся в норме
 пространства Соболева $\mathbf{H}^s(B)$, необходимо и достаточно,
  чтобы $\mathbf{f}$ принадлежала
 \[\mathbf{V}^s_\mathcal{R}(B)=
\{\mathbf{f}\in \mathbf{V}^0\cap\mathbf{H}^s(B):
\gamma_{\mathbf{n}}\mathbf{f}=0,...,
 \gamma_{\mathbf{n}} \mathrm{rot}^{s-1}\mathbf{f}=0\},\,\,
 \|\mathbf{f}\|_{\mathbf{V}_\mathcal{R}^s}=
 \|\mathbf{f}\|_{\mathbf{H}^s}.\]

 Если $\mathbf{f}\in \mathbf{V}^s_\mathcal{R}(B)$,
то сходится ряд
\begin{equation} \label{rof 2}
\sum_{\kappa, n>0}{\lambda}_{\kappa}^{2s}\,
(|(\mathbf{f},\mathbf{q}_{\kappa}^+)|^2
+|(\mathbf{f},\mathbf{q}_{\kappa}^-|^2)),\quad
{\lambda}_{\kappa}=({\rho}_{n,m})/R
\end{equation} и существует такая , не
зависящая от $\mathbf{f}$ постоянная $C>0$, что
\begin{equation} \label{orf 3}
\sum_{\kappa, n>0} {\lambda}_{\kappa}^{2s}\,
(|(\mathbf{f},\mathbf{q}_{\kappa}^+)|^2
+|(\mathbf{f},\mathbf{q}_{\kappa}^-|^2))\leq
C\|\mathbf{f}\|^2_{\mathbf{H}^s(B)}.
\end{equation}

 Если $s\geq 2$, то любая вектор-функция  $\mathbf{f}$ из
  $\mathbf{V}^s_\mathcal{R}(B)$
разлагается в в ряд Фурье, сходящийся в пространстве $\mathbf{C}^{s-2}(\overline{B})$.%
\end{theorem}

\section{Базисные подпространства в $\mathbf{L}_{2}(G)$}

\subsection{Подпространства 
${{\mathbf{L}}_{2}}({G})$} Как мы показали в $\S 1$ пространство
$\mathbf{L}_{2}(G)$ допускает разложение \eqref{bhr  1}
  на ортогональные  подпространства $\mathcal{A}$, 
$\mathcal{B}_H$ и $\mathbf{V}^0$, которое запишем так:
\begin{equation}
\label{wro 1}\mathbf{L}_{2}(G)=\mathcal{A}\oplus\mathcal{B}_H\oplus
\mathbf{V}^{0}(G).
\end{equation}
Пространства $\mathcal{A}$ и $\mathcal{B}_H$ принадлежат ядру
оператора ротор; а в $\mathbf{V}^0$ он продолжается как
самосопряженный оператор $S$, собственные функции
$\mathbf{q}^{\pm}_{k }(\mathbf{x})$ которого образуют
ортонормированный в $\mathbf{V}^0$ базис.

Согласно $\S 2$  пространство $\mathcal{B}=\mathcal{B}_H\oplus\
\mathbf{V}^0$ принадлежит ядру оператора градиент дивергенции; а в
$\mathcal{A}_\gamma\subset \mathcal{A}$ он продолжается как
самосопряженный оператор $\mathcal{N}_d$, собственные функции
$\mathbf{q}_{j }(\mathbf{x})$ которого образуют ортонормированный в
$\mathcal{A}$ базис.

 $\mathcal{B}_H$--конечномерное пространство,
  $\mathbf{h}_{i }(\mathbf{x})$--его базисные   функции.

Следовательно, имеет место
\begin{theorem} Система векторов
$\{\mathbf{q}_{j}(\mathbf{x})\}\cup\{\mathbf{h}_i(\mathbf{x})\}\cup\{\mathbf{q}_{k
}^{+}(\mathbf{x})\}\cup \{\mathbf{q}_{k }^{-}(\mathbf{x})\}$ 
образует в пространстве ${\mathbf{{L}}_{2}}(G)$ ортонормированный
базис. Любую вектор-функцию  $\mathbf{f}(\mathbf{\mathbf{x}})$ из
${\mathbf{{L}}_{2}}(G)$
 можно разложить в ряд Фурье по этому базису:
\[\mathbf{f}(\mathbf{x})=\sum_{\mu_j\in M}(\mathbf{f},\mathbf{q}_{j})\mathbf{q}_{j}
(x)+ \sum_{i=1}^{\rho}(\mathbf{f},\mathbf{h}_{i})
 \mathbf{h}_{i}(\mathbf{x})+\sum_{\pm\lambda_k\in
\Lambda}[(\mathbf{f},\mathbf{q}^{+}_{k})\mathbf{q}^{+}_{k}+
(\mathbf{f},\mathbf{q}^{-}_{k})\mathbf{q}^{-}_{k}].
\]
 \end{theorem}

 В случае шара пространство $\mathcal{B}_H$ пусто.
\subsection{Разложение векторного поля $\mathbf{f}$ из $\mathbf{L}_{2}(B)$
 на потенциальное
 и соленоидальное поля $\mathbf{f}_{\mathcal{A}}$ и $\mathbf{f}_{\mathcal{B}}$,}\quad
$\mathbf{f}(\mathbf{x})=\mathbf{f}_{\mathcal{A}}(\mathbf{x})+\mathbf{f}_{\mathcal{B}}(\mathbf{x})$,
где
\begin{equation} \label{sra 1}
\mathbf{f}_{\mathcal{A}}=
\sum_{n=0}^{\infty}\sum_{m=1}^{\infty}\sum_{k=-n}^{n}
(\mathbf{f},{\mathbf{q}}_{\kappa})\,\mathbf{q}_{\kappa}(\mathbf{x}),\quad
\kappa=(n,m,k)\end{equation}

\begin{equation} \label{srpm 2}
\mathbf{f}_{\mathcal{B}}=\sum_{n=1}^{\infty}\sum_{m=1}^{\infty}\sum_{k=-n}^{n}
[(\mathbf{f},{\mathbf{q}}_{\kappa}^+)\,\mathbf{q}_{\kappa}^+(\mathbf{x})
+(\mathbf{f},{\mathbf{q}}_{\kappa}^-)\,\mathbf{q}_{\kappa}^-(\mathbf{x})].
\end{equation}
{\small Частичные суммы $\mathbf{S}^0_N$ и $\mathbf{S}^1_N$ рядов
\eqref{sra 1} и \eqref{srpm 2} состоят из элементов с индексами
$\kappa$, для которых $0<\alpha_{n,m}<N$ и $0<\rho_{n,m}<~N$,
соответственно. Числа $\rho_{n,m}$ и $\alpha_{n,m}$ - нули функций
$\psi_{n}(x)$ и их производных.}

 При $\mathbf{f},\mathbf{g}\in \mathbf{L}_{2}(B)$
 \, скалярное произведение \[ (\mathbf{f}, \mathbf{g})=
 ( \mathbf{f}_{\mathcal{A}},\mathbf{g}_{\mathcal{A}})+
 ( \mathbf{f}_{\mathcal{B}},\mathbf{g}_{\mathcal{B}}).\] Равенство Парсеваля-Стеклова,
$\|\mathbf{f}\|^2=\|\mathbf{f}_{\mathcal{A}}\|^2+\|\mathbf{f}_{\mathcal{B}}\|^2$,
 запишем так
\begin{equation} \label{psra 1}
\|{\mathbf{f}}\|^2=\sum_{\kappa, n\geq 0}
(\mathbf{f},{\mathbf{q}}_{\kappa})^2 +\sum_{\kappa,n\geq
1}[(\mathbf{f},{\mathbf{q}}_{\kappa}^+)^2
+(\mathbf{f},{\mathbf{q}}_{\kappa}^-)^2],
\end{equation}

Отметим, что разложение векторного поля $\mathbf{f}(\mathbf{x})$ на
потенциальное поле $\nabla{h}(\mathbf{x})$ и соленоидальное поле
$\mathbf{u}(\mathbf{x})$ связано  с решением задачи Неймана
\begin{equation} \label{npr 1}
\triangle h= \text{div}\,\mathbf{f} \quad \text{в}\quad B,\quad
\mathbf{n}\cdot\nabla\,h|_S=\mathbf{n}\cdot\mathbf{f}|_S,
\end{equation}в классической или обобщенной постановках
(cм.\cite{lad, bs}).
 Мы  находим решение  задачи
 в виде рядов \eqref{sra 1}, \eqref{srpm 2}.
\subsection {Свойства рядов \eqref{sra 1}, \eqref{srpm 2}} Пространство
$\mathbf{L}_{2}(B)\subset \mathbf{L}_{1}(B)$, а скалярное
произведение $ (\mathbf{f}, \mathbf{g})$ есть регулярный функционал
в  $\mathcal{D}'(B)$. Ряды  $\mathbf{f}_{\mathcal{A}}$ и
$\mathbf{f}_{\mathcal{B}}$  задают распределения
$(\mathbf{f}_{\mathcal{A}}, \,\mathbf{g}_{\mathcal{A}})$ и
$(\mathbf{f}_{\mathcal{B}}, \,\mathbf{g}_{\mathcal{B}})$ при
$\mathbf{g}_{\mathcal{A}}$ и $\mathbf{g}_{\mathcal{B}}$ из
$\mathcal{D}(B)$. Их дифференцируют почленно. Частичные суммы
$\mathbf{S}^0_N$ и $\mathbf{S}^1_N$  удовлетворяют уравнениям:
$\mathrm{rot}\,\mathbf{S}^0_N=0,\,
\mathrm{div}\,\mathbf{S}^1_N=~0$\, и условиям \,
$\gamma_{\mathbf{n}}\,\mathbf{S}^0_N=0,\,
\gamma_{\mathbf{n}}\,\mathbf{S}^1_N=0.$
 Поэтому
$\mathrm{rot}\,\mathbf{f}_{\mathcal{A}}=0$ \,и \,
$\mathrm{div}\,\mathbf{f}_{\mathcal{B}}=0$ в $\mathcal{D}'(B)$.

  Согласно Теореме 3 п.\,3.6 ряд \eqref{sra 1} сходится в
  $\mathbf{H}^{1}(B)$ (соотв., в $\mathbf{H}^{2}(B)$), если
  $\mathbf{f}\in\mathbf{H}^{1}(B)$ (соотв.,
    $\mathbf{f}\in\mathbf{H}^{2}(B)$) и
    $\gamma_{\mathbf{n}}\mathbf{f}=0$.

    Ряд сходится в
     $\mathbf{H}^{3}(B)$, а, значит, и в $\mathcal{C}^{1}(\overline{B})$,
      если
     $\mathbf{f}\in\mathbf{H}^{3}(B)$ и
     $\gamma_{\mathbf{n}}\mathbf{f}=0$,
      $\gamma_{\mathbf{n}}\nabla\mathrm{div}\mathbf{f}=0$.

      Обратно, если ряд \eqref{sra 1} сходится в
  $\mathbf{H}^{1}(B)$ (соотв., в $\mathbf{H}^{2}(B)$), то
   $\gamma_{\mathbf{n}}\mathbf{f}_{\mathcal{A}}=0$.
  Если ряд $\mathbf{f}_{\mathcal{A}}$ сходится в
  $\mathbf{H}^{3}(B)$ (соотв., в    $\mathbf{H}^{4}(B)$), то\newline
   $\gamma_{\mathbf{n}}\mathbf{f}_{\mathcal{A}}=0$ и
   $\gamma_{\mathbf{n}}\nabla\mathrm{div}\mathbf{f}_{\mathcal{A}}=0$.
    И так далее.

\subsection {Пространства $ \mathbf{E}^{s}(B)$, $ \mathbf{F}^{s}(B)$}
Через  $\mathbf{H}^{s}_{div}(B)$\,
 или $
\mathbf{E}^{s}(B)$ \,\cite{rt}
обозначают пространство   \,\,   
\[\mathbf{E}^{s}(B)=\{\mathbf{v}\in
\mathbf{H}^{s}(B):\,\text{div}\mathbf{v}\in {H}^{s}(B)\}\]\,\,\,
 с нормой\,\,$\|\mathbf{v}\|_{\mathbf{E}^{s}}=(
 \|\mathbf{v}\|^2_{s}+
 \|\text{div}\,\mathbf{v}\|^2_{s})^{1/2},$
  где число
 $s\geq 0$ целое.

 Оно является пространством Гильберта и
\[\mathbf{C}_0^{\infty}({B})  \subset \mathbf{E}^{s}(B), \quad
 \mathbf{H}^{s+1}(B)\subset \mathbf{E}^{s}(B)
 \subset\mathbf{H}^s(B). \]

В  терминах проекций
$\mathbf{f}_{\mathcal{A}},\,\mathbf{f}_{\mathcal{B}}$ условие
$\mathbf{f}\in \mathbf{E}^{0}(B)$ означает, что
\[\mathbf{f}_{\mathcal{A}},\,\mathbf{f}_{\mathcal{B}},\,
\mathrm{div}\,\mathbf{f}_{\mathcal{A}}\in \mathbf{L}_{2}(B)\quad
\text{и} \quad\|\mathbf{f}\|_{\mathbf{E}^{0}}^2=
 \|\mathbf{f}\|^2+
 \|\text{div}\,{\mathbf{f}_{\mathcal{A}}}\|^2,\]
 \[\|\text{div}\,\mathbf{f}_{\mathcal{A}}\|^2= \sum_{N=1}^{\infty} \sum_{(n,m)\in
\mathbb{P}^1_N}^{}\sum_{k\in[-n,n]}^{}(\alpha_{n,m}/R)^2
(\mathbf{f},{\mathbf{q}}_{\kappa})^2 ,\]где
$\mathbb{P}^1_N=\{(n,m):0<\alpha_{n,m}<N\}$,\quad
$\psi^{'}_n(\alpha_{n,m})=0$.

Очевидно, что ${\bf rot}\, \mathbf{u}+\lambda \mathbf{u}\in
\mathbf{E}^{s}(B)$, если $\mathbf{u}\in \mathbf{H}^{s+1}(B)$.

 Для функций $v$ из пространства $H^1(B)$ определен \,\cite{mi}
оператор {\it следа} $\gamma :H^1(B)\rightarrow H^{1/2}(S)$, равный
следу $v$ на $S$ для гладких функций из
$\mathcal{C}^1(\overline{B})$: $\gamma\,v=v|_S$,
 причем $\quad \|\gamma\,v\|_{L_2(S)}\leq c \|v\|_{H^1(B)}$.

 Аналогично, для вектор-функций $\mathbf{u}(\mathbf{x})$ из
$\mathbf{E}^0(B)$ определен \cite{rt} оператор {\it следа нормальной
компоненты} $\gamma_\mathbf{n} :\mathbf{E}^0(B)\rightarrow
H^{-1/2}(S)$, равный сужению $\mathbf{n}\cdot \mathbf{u}$ на $S$ для
функций из $\mathcal{C}(\overline{B})$:
$\gamma_\mathbf{n}\,\mathbf{u}=\mathbf{n}\cdot\mathbf{u}|_S$.

 Для $\mathbf{u}\in \mathbf{E}^0(B)$ и $v\in H^1(B)$ верна
обобщенная формула Стокса:
 $\langle \gamma_\mathbf{n}{\mathbf{u}},\gamma\,v \rangle=
 (\mathbf{u}, \nabla v)+( \text{div}\mathbf{u}, v),$
  где  $\langle\gamma_\mathbf{n}{\mathbf{u}},\gamma\,v \rangle$-
  линейный
функционал над $H^{1/2}(S)$. Имеют место 
вложения:
 $H^{1/2}(S)\subset
L_2(S)\subset H^{-1/2}(S).$

 Аналогично, согласно п.2.2
\[{\mathbf{F}^{s}}(B)=\{\mathbf{v}\in {\mathbf{H}^{s}(B)}:
\text{rot}^2\,\mathbf{v}\in
 \mathbf{H}^{s}(B)\}\] это пространство с нормой
 $\|\mathbf{v}\|_{\mathbf{F}^{s}}=( \|\mathbf{v}\|^2_{s}+
 \|\mathrm{curl}^2\mathbf{v}\|^2_{s})^{1/2}$.

В  терминах проекций
$\mathbf{f}_{\mathcal{A}},\,\mathbf{f}_{\mathcal{B}}$ условие
$\mathbf{f}\in \mathbf{F}^{0}(B)$ означает, что
\[\mathbf{f}_{\mathcal{A}},\,\mathbf{f}_{\mathcal{B}},\,
\mathrm{rot}^2\,\mathbf{f}_{\mathcal{B}}\in \mathbf{L}_{2}(B)\quad
\text{и} \quad\|\mathbf{f}\|_{\mathbf{F}^{0}}^2=
 \|\mathbf{f}\|^2+ 
 \|\text{rot}^2\,{\mathbf{f}_{\mathcal{B}}}\|^2,\]
 \[\|\text{rot}^2\,\mathbf{f}_{\mathcal{B}}\|^2= \sum_{N=1}^{\infty} \sum_{(n,m)\in
\mathbb{P}^{2}_N}^{}\sum_{k\in[-n,n]}^{}(\rho_{n,m}/R)^4
[(\mathbf{f},{\mathbf{q}}_{\kappa}^+)^2
+(\mathbf{f},{\mathbf{q}}_{\kappa}^-)^2 ],\]где
$\mathbb{P}^2_N=\{(n,m):0<\rho_{n,m}<N\}$,\quad
$\psi^{}_n(\rho_{n,m})=0$\,\, и \,\,${\mathbf{q}}_{0,m,0}^{\pm}=0$.

Пространства с индексом ${\gamma}$  \,\, определяются так:
\[\mathbf H^1_{\gamma}(B)=
\{\mathbf{f}\in\mathbf{H}^l(B): \gamma_\mathbf{n}\mathbf{f}=0\},
\quad \mathbf{H}^2_{\gamma\delta\gamma}(B)=\{\mathbf{g}\in
\mathbf{{H}}^2 _{\gamma}(B):\,\gamma_
 {\mathbf{n}}{\nabla\text{div}\,\mathbf{g}}=0 \},\]
\[\mathbf{{E}}_{\gamma}^s(B)=
 \{\mathbf{f}\in\mathbf{E}^s(B):
\gamma_\mathbf{n}{\mathbf{f}}=0 \}, \quad s\geq 0,\quad \mathbf
H^{s+1}_{\gamma}\subset\mathbf E^s_{\gamma}\subset H^{s}_{\gamma},\]
 \[\mathbf{{F}}_{\gamma}^s(B)=
 \{\mathbf{f}\in\mathbf{F}^s(B):
\gamma_\mathbf{n}{\mathbf{f}}=0 \}, \quad s\geq 0,\quad \mathbf
H^{s+2}_{\gamma}\subset\mathbf
F^s_{\gamma}\subset\,H^{s}_{\gamma}.\]

\section{Решение  краевой задачи  в шаре}

 Методом Фурье легко
решается  краевая задача.

\begin{problem} Задана вектор-функция $\mathbf{f}(\mathbf{x})\in
\mathbf{L}_{2}(B)$. Найти вектор-функцию $\mathbf{v}(\mathbf{x})$ в
$\mathbf{H}^{2}(B)$  такую, что
 \begin{equation}
\label{grad__3_} {\nabla\bf div}\, \mathbf{v}+\lambda \mathbf{v}=
\mathbf{f}\quad \text{в}\quad B, \quad \gamma_\mathbf{n}
\mathbf{v}=0,
\end{equation}\end{problem}
Определение.   {\it Вектор-функция $\mathbf{v}$ из
$\mathbf{L}_{2}(B)$ есть обобщенное решение задачи при
$\mathbf{f}\in\mathbf{L}_{2}(B)$, если она удовлетворяет тождеству }
 \begin{equation}
\label{ogr__2_}(\mathbf{v}, (\nabla\,\mathrm{div}\, +\lambda
)\mathbf{w})=
 (\mathbf{f},\mathbf{w})\quad \text{для \, любой }\quad
 \mathbf{w}\in \mathbf{H}^{2}_0(B).\end{equation}

 Если $\mathbf{f}=\mathbf{f}_{\mathcal{B}}\in \mathcal{B}$
и $\lambda\neq 0$, то
$\mathbf{v}=~{\lambda}^{-1}\mathbf{f}_{\mathcal{B}}$
   есть решение уравнения  \eqref{ogr__2_}.


 \subsection {Решение  краевой задачи 3
при  $\lambda\,\overline{\in}\, Sp\, (-\nabla\mathrm{div})$}  Имеет
место
\begin{theorem}Если $\lambda\neq 0, {\nu}_{n,m}^2$ $(n\geq 0,\,m>0)$,\,
$\mathbf{f}(\mathbf{x})\in {\mathbf{F}^0_{\gamma}}(B)$,\,
  то единственное решение   $\mathbf{v}$ задачи 3 есть сумма
  $\mathbf{v}_1+\mathbf{v}_2$ рядов
 \begin{equation} \label{gra 1}
{\mathbf{v}_1}=\sum_{n=0}^{\infty}\sum_{m=1}^{\infty}\sum_{k=-n}^{n}
\frac{(\mathbf{f},{\mathbf{q}}_{\kappa})}{\lambda-\nu_{n,m}^2}\,
\mathbf{q}_{\kappa} (\mathbf{x}),\quad \kappa=(n,m,k),
\end{equation}
\begin{equation} \label{gra 2}\mathbf{v}_2=
\frac{\mathbf{f}_{\mathcal{B}}}{\lambda}\equiv
\frac{1}{\lambda}\sum_{n=1}^{\infty}\sum_{m=1}^{\infty}\sum_{k=-n}^{n}
[(\mathbf{f},{\mathbf{q}}_{\kappa}^+)\,\mathbf{q}_{\kappa}^+
+(\mathbf{f},{\mathbf{q}}_{\kappa}^-)\,
\mathbf{q}_{\kappa}^-(\mathbf{x})].\end{equation} Решение задачи
принадлежит пространству Соболева $\mathbf{H}^{2}_\gamma(B)$.

Если   $\mathbf{f}\in \mathcal{A}\subset\, \mathbf{L}_{2}(B)$, то
 $\mathbf{v}=(\mathcal{N}_{d}\,+\lambda)^{-1}\,\mathbf{f}$ принадлежит
$\mathcal{A}^{2}_{\gamma}\subset
\mathbf{H}^{2}_{\gamma}$.

Если $\mathbf{f}\in \mathcal{D}(B)$, то ряды \eqref{gra 1} и
\eqref{gra 2} Сходятся в любом из пространств
 $\mathbf{H}^{s}(B)$,\, $s\geq 1$, и их сумма есть
 решение  задачи класса $C^{\infty}(\overline{B})$.
\end{theorem}
{\it 1. Существование и единственность обобщенного решения.}

 Пусть  $\mathbf{v}\in
\mathbf{H}^{2}_\gamma(B)$ и
  $\mathbf{f}\in \mathbf{F}^{0}_\gamma(B)$ в \eqref{grad__3_}.\,
  Уравнение \eqref{grad__3_}
умножим на $\mathbf{w }\in\mathbf{H}^{2}_0(B)$ и проинтегрируем по
частям. Получим равенство \eqref{ogr__2_}.

 Откуда,
 полагая $\mathbf{w}=\mathbf{q}_{\kappa}
(\mathbf{x})$ и  $\mathbf{w}=\mathbf{q}_{\kappa}^{\pm}
(\mathbf{x})$, соответственно,  получим соотношения:
\begin{equation} \label{fogra 2}
(\lambda-\nu_{n,m}^2)\, (\mathbf{v},{\mathbf{q}}_{\kappa})=
(\mathbf{f},{\mathbf{q}}_{\kappa}),\quad  \lambda \,
(\mathbf{v},{\mathbf{q}}_{\kappa}^{\pm})=
(\mathbf{f},{\mathbf{q}}_{\kappa}^{\pm}).\end{equation}

 Формулы \eqref{gra 1} и \eqref{gra 2}
  получают  из  соотношений \eqref{fogra 2}.

Ряды \eqref{gra 1},\,\eqref{gra 2} сходятся  в $\mathbf{L}_{2}(B)$,
так как по условию $\mathbf{f}\in\mathbf{L}_{2}(B)$ и числа
$|\lambda|^{-1}$, $|\lambda-\nu_{n,m}^2)|^{-1}$ ограничены
    при  любых $n,m$.

Пусть $\mathbf{w}\in \mathcal{D}(B)$. Функция
$\mathbf{v}_2=~{\lambda}^{-1}\mathbf{f}_{\mathcal{B}}$ есть решение
уравнения \eqref{ogr__2_}: $({\lambda}^{-1}\mathbf{f}_{\mathcal{B}},
(\nabla\mathrm{div}\,\mathbf{w}+\lambda\,\mathbf{w}))=
(\mathbf{f}_{\mathcal{B}},\mathbf{w})=
(\mathbf{f}_{\mathcal{B}},\mathbf{w}_{\mathcal{B}})$. При
$\mathbf{w}=\mathbf{w}_{\mathcal{A}}$
\[\nabla\mathrm{div}\,\mathbf{w}+\lambda\,\mathbf{w}=
\sum_{\kappa,n\geq
0}(\mathbf{w},{\mathbf{q}}_{\kappa})(\lambda-\nu_{n,m}^2)\,
\mathbf{q}_{\kappa}(\mathbf{x}) ,\]
\[(\mathbf{v}_{1},
(\nabla\mathrm{div}\,\mathbf{w}+\lambda\,\mathbf{w}))
=\sum_{\kappa,n\geq 0}[(\mathbf{f},{\mathbf{q}}_{\kappa})
(\mathbf{w},{\mathbf{q}}_{\kappa})=
(\mathbf{f}_{\mathcal{A}},\mathbf{w}_{\mathcal{A}} ) .\]

Следовательно, $((\mathbf{v}_{1}+\mathbf{v}_{2}),
(\nabla\mathrm{div}\,\mathbf{w}+\lambda\,\mathbf{w}))=
(\mathbf{f},\mathbf{w})$.

Существование обобщенного решения доказано.

Единственность решения задачи 3  вытекает из
 полноты семейства собственных функций ротора  и градиента дивергенции
 в $\mathbf{L}_2(B)$.

 {\it 2. Операторы задачи и сходимость рядов.}  Обозначим ряд
 \eqref{gra 1}  через $Q_{\lambda}^{}\,{\mathbf{f}}$, а
  ряд \eqref{gra 2}--
  через ${\lambda}^{-1}\mathcal{P}_\mathcal{B}\,{\mathbf{f}}$.

Спектр $ \{0,  {\nu}_{n,m}^2\}$ оператора $-\nabla\mathrm{div}$ не
имеет конечных предельных точек, поэтому числа\,$|\lambda |^{-1}$,\,
$|\lambda -{\nu}_{n,m}^2|^{-1}$ ограничены и
\begin{equation} \label{ner 2}
\|{\lambda}^{-1}\mathcal{P}_\mathcal{B}\,{\mathbf{f}}\|\leq
|{\lambda}|^{-1}\|{\mathbf{f}}\|, \quad
\|Q_{\lambda}^{}\,{\mathbf{f}}\| \leq \Lambda
 \|{\mathbf{f}}\|,  \quad \Lambda= \max_{n,m}\,1/
 |{\lambda}- {\nu_{n,m}^2|}
\end{equation}
причем постоянные $|\lambda|^{-1}$ и $\Lambda$ зависят только от
расстояния точки $\lambda$ от точек спектра градиента дивергенции.

 Пусть $\mathbf{f}$ принадлежит пространству
 $\mathbf{F}_{\gamma}^0(B)$, то-есть $\mathbf{f}$ и
 $\text{rot}^2\,\mathbf{f}\in L_{2}(B)$.
 Тогда $\mathbf{v}_2={\lambda}^{-1}\mathcal{P}_\mathcal{B}\,
 {\mathbf{f}}$ и
$\text{rot}^2\,\mathbf{v}_2={\lambda}^{-1}\text{rot}^2\,
\mathbf{f}_\mathcal{B}= {\lambda}^{-1}\text{rot}^2\,\mathbf{f}$\,
также принадлежат $L_{2}(B)$.
 Так как $\text{div}\,\mathbf{v}_2=0$ в $B$,
 \, $\gamma_{\mathbf{n}}\,\mathbf{v}_2=0$, то
 согласно оценкам  \eqref{ner 2} и \eqref{apgro s}\,(при $s=0$)
\begin{equation} \label{apr 2}
\|{\lambda}^{-1}\mathcal{P}_\mathcal{B}\,{\mathbf{f}}\|_2\leq
(C_0|{\lambda}|)^{-1}(\|{\mathbf{f}}\|+\|\text{rot}^2\,\mathbf{f}\|).
\end{equation}
Значит ${\lambda}^{-1}\mathcal{P}_\mathcal{B}\,{\mathbf{f}}\in
\mathbf{H}^2_{\gamma\delta\gamma}(B)$ и
$(\nabla\text{div}+\lambda)\,{\lambda}^{-1}\mathbf{f}_\mathcal{B}=
 \mathbf{f}_\mathcal{B}$.

Далее, вектор-функция $\mathbf{v}_1(\mathbf{x})=
Q_{\lambda}\,\mathbf{f}$ принадлежит  пространству
\[\mathcal{A}_\gamma(B)=\{\mathbf{v}\in \mathcal{A}:
\gamma_{\mathbf{n}}\mathbf{v} =0\}.\]
  \begin{equation} \label{grad 2}
\nabla \text{div}\, \mathbf{v}_1
=-\sum_{n=0}^{\infty}\sum_{m=1}^{\infty}\sum_{k=-n}^{n}
\frac{{\nu}^2_{n,m}}{\lambda - {\nu}^2_{n,m}} (\mathbf{f},
{\mathbf{q}}_{\kappa})\, {\mathbf{q}}_{\kappa}(\mathbf{x}),\quad
\kappa=(n,m,k),
\end{equation}
\begin{equation} \label{pi q}
 \|\nabla \text{div}\, \mathbf{v}_1
\|_{}^2
\leq \Pi^2\,\|{\mathbf{f}}\|^2,
 \,\,\text{где}\,\,\Pi^2= \max_{n,m}\, \frac{\nu_{n,m}^4}{|\lambda -
\nu^2_{n,m}|^{2}}. \end{equation} Следовательно,
 $\mathbf{v}_1\in\mathcal{A}^2_\gamma(B)$ и из оценок \eqref{ner 2}
 и \eqref{pi q}
вытекает, что
\begin{equation} \label{apr q}
\|Q_{\lambda}\,{\mathbf{f}}\|_2\leq
C_0^{-1}(\Lambda+\Pi)\|{\mathbf{f}}\|.
\end{equation}
Значит, $Q_{\lambda}\,{\mathbf{f}}\in
 \mathbf{{H}}^2_{\gamma\delta\gamma}(B)$ и
 $(\nabla\text{div}+\lambda)\,Q_{\lambda}\,{\mathbf{f}}=
 \mathbf{f}_\mathcal{A}$.

  Поэтому
  $\mathbf{v}(\mathbf{x})=
 {\lambda}^{-1}\mathcal{P}_\mathcal{B}\,{\mathbf{f}}+
 Q_{\lambda}\,{\mathbf{f}}$ есть
 решение задачи 3.
 Оно принадлежит пространству
 $\mathbf{{H}}^2_{\gamma\delta\gamma}(B)$ и
\begin{equation} \label{apr 7}
\|v\|_{\mathbf{{H}}^2_{\gamma\delta\gamma}(B)}\leq
C_0^2\|{\mathbf{f}}\|_{\mathbf{{F}}_{\gamma}^0(B)}.
\end{equation}
Если $\mathbf{f}\in\mathbf{F}_{\gamma}^{s}(B)$, $s\geq\,3$, то ряды
\eqref{gra 1}, \eqref{gra 2} сходятся в  пространстве
$C^{s-1}(\overline{B})$ и их сумма есть классическое решение задачи
класса $C^{s-1}(\overline{B})$.
 Наконец, если $\mathbf{f}\in \mathcal{D}(B)$, то согласно п.4.3
ряды \eqref{gra 1}, \eqref{gra 2} сходятся в любом из пространств
$\mathbf{H}^{s}(B)$,\, $s\geq1,$ и задают 
решение задачи класса $C^{\infty}(\overline{B})$.
   Теорема доказана.
Кроме того, имеет место
\begin{lemma}   Если $\lambda\,\overline{\in}\, Sp\,
(-\nabla\mathrm{div})$,
то оператор $\nabla\text{div}+\lambda I$ осуществляет гомеоморфизм
(т.-е. взаимно однозначное и непрерывное
 отображение) пространств  $\mathbf{{H}}^2_{\gamma\delta\gamma}(B)$
 и  $\mathbf{{F}}_{\gamma}^0(B)$.\end{lemma}
Действительно, пусть
$\mathbf{v}\in\mathbf{{H}}^2_{\gamma\delta\gamma}
=\{\mathbf{v}\in
\mathbf{{H}}^2:\,\gamma_{\mathbf{n}}\,\mathbf{v}=0,
\gamma_{\mathbf{n}}\nabla\mathrm{div}\,\mathbf{v}=0\}$, тогда
$\mathbf{f}=\nabla\,\text{div}\,\mathbf{v}+\lambda\mathbf{v}\in
\mathbf{L}_{2}(B)$,\,
$\text{rot}^2\mathbf{f}={\lambda}\text{rot}^2\mathbf{v}\in
{L}_{2}(B)$ и \newline $\gamma_{\mathbf{n}}\mathbf{f}=
\gamma_{\mathbf{n}}\nabla\text{div}\mathbf{v}+
\lambda\gamma_{\mathbf{n}}\mathbf{v}=0$.
 Значит, ${\mathbf{f}}\in \mathbf{{F}}_{\gamma}^0(B)$
 и 
\begin{equation} \label{opr 6}
\|f\|_{\mathbf{{F}}_{\gamma}^0(B)}\leq
C_1^0\|{\mathbf{v}}\|_{\mathbf{{H}}^2_{\gamma\delta\gamma}(B)}.
\end{equation}
Обратный оператор
${\lambda}^{-1}\mathcal{P}_\mathcal{B}\,{\mathbf{f}}+
 Q_{\lambda}^{}\,{\mathbf{f}}$ отображает
 $\mathbf{{F}}_{\gamma}^0(B)$ на
 $\mathbf{{H}}^2_{\gamma\delta\gamma}(B)$ и 
 согласно оценке \eqref{apr 7} непрерывен. Лемма 
 доказана

\subsection {Разрешимость  краевой задачи 3
при  $\lambda\in Sp\, (-\nabla\mathrm{div})$} Из соотношений
\eqref{fogra 2} видим, что

{\it при $\lambda=0$ однородная задача имеет счетное
 линейно независимых решений $\mathbf{q}_{\kappa}^{\pm}(\mathbf{x})$,
а неоднородная задача 5 разрешима тогда и только тогда, когда
 $(\mathbf{f},{\mathbf{q}}_{\kappa}^{\pm})=0\,\, \forall\, \kappa$, то-есть
 $\mathbf{f}_{\mathcal{B}}=0$    или  $\text{rot}\,\mathbf{f}=0$,

при $\lambda=\nu_{n,m}^2$ и фиксированных $n, m$ однородная задача
имеет $2n+1$
 линейно независимых решений $\mathbf{q}_{\kappa}(\mathbf{x})$, где
$\kappa=(n,m,k),\,\, k=-n,...,n$, а неоднородная задача 5 разрешима
тогда и только тогда, когда
 $(\mathbf{f},{\mathbf{q}}_{\kappa})=0$, то-есть задача разрешима по Фредгольму.}

\subsection{Оператор $\mathcal{N}_{d}+\lambda\, I$ в
 подпространстве $\mathcal{A}_{\gamma}$ }
Согласно п.2.3 оператор $\nabla\mathbf{div}+\lambda\,I$ и его
расширение $\mathcal{N}_{d}+\lambda\,I$ в $\mathcal{A}_{\gamma}$
задаются рядом
\[(\mathcal{N}_{d}\,+\lambda)\,\mathbf{w}=
\sum_{\kappa,n\geq
0}(\mathbf{w},{\mathbf{q}}_{\kappa})(\lambda-\nu_{n,m}^2)\,
\mathbf{q}_{\kappa}(\mathbf{x}),\quad \mathbf{w}\in
\mathcal{A}_{\gamma} ,\] если он сходится в $\mathbf{L}_2(B)$.
Обратный оператор имеет вид:
\[(\mathcal{N}_{d}\,+\lambda)^{-1}\,\mathbf{v}=
\sum_{\kappa,n\geq
0}\frac{(\mathbf{v},{\mathbf{q}}_{\kappa})}{(\lambda-\nu_{n,m}^2)}\,
\mathbf{q}_{\kappa}(\mathbf{x}),\quad \mathbf{v}\in
\mathcal{A}_{\gamma},\] если $\lambda\neq\nu_{n,m}^2$, где
  $\nu_{n,m}=({\alpha }_{n,m})/R$,    $n\ge 0$,   $m\in N$, а   числа
      ${{\alpha }_{n,m}}>0$ суть нули  функций
      ${{\psi }_{n}^{\prime }}(z)$,
     производных функций \eqref{psi__1_}. 


Если $\lambda=\nu_{n,m}^2$ при фиксированных $n$  и $m$, то
однородное уравнение $(\mathcal{N}_{d}\,+\lambda)\,\mathbf{w}= 0$
имеет $2n+1$
 линейно независимых решений $\mathbf{q}_{\kappa}(\mathbf{x})$, где
$\kappa=(n,m,k),\,\, k=-n,...,n$, которые являются собственными
функциями оператора  $-\nabla\mathbf{div}$ и вычислены явно в п.3.3.
Неоднородное уравнение $(\mathcal{N}_{d}\,+\lambda)\,\mathbf{w}=
\mathbf{v}$ разрешимо тогда и только тогда, когда
 $(\mathbf{v},{\mathbf{q}}_{\kappa})=0$ \,при\,
 $\kappa=(n,m,k),\,\, k=-n,...,n$.

 Значит оператор
 $\mathcal{N}_{d}\,+\lambda\,I:
 \mathcal{A}_{\gamma}\rightarrow\mathcal{A}_{\gamma}$
 является  фредгольмовым, а оператор
 $\mathcal{N}_{d}:
 \mathcal{A}_{\gamma}^2\rightarrow\mathcal{A}_{\gamma}$ --
 однозначно обратимым.

Определены степени операторов $\mathcal{N}_{d}$ и
$\mathcal{N}_{d}^{-1}$:
\[\mathcal{N}_{d}^{\pm p}\,\mathbf{v}= \sum_{\kappa,n\geq
0}(\mathbf{v},{\mathbf{q}}_{\kappa})(-\nu^2_{n,m})^{\pm p}\,
\mathbf{q}_{\kappa}(\mathbf{x}),\quad p=2,3,....\]
 По Теореме 3 ряды $\mathcal{N}_{d}^{p}$ сходятся в $\mathbf{L}_2(B)$ тогда и только
тогда, когда
\[\mathbf{v}\in\mathbf{A}^p_\mathcal{K}(B)=\{\mathbf{f}\in
\mathcal{A}\cap\mathbf{H}^p(B): \gamma_{\mathbf{n}}\mathbf{f}=0,...,
 \gamma_{\mathbf{n}}
  (\nabla\mathrm{div})^{\widehat{p}}\mathbf{f}=0\},\,\,\widehat{p}=
  [p/2].\]
 Степени обратного оператора $\mathcal{N}_{d}^{-1}$
отображают пространство $\mathcal{A}_{\gamma}(B)$ на пространства
$\mathcal{A}_{\gamma}^p(B)=\{\mathbf{v}\in
\mathcal{A}_{\gamma}(B),...,(\nabla\mathrm{div})^p\mathbf{v}\in
\mathcal{A}_{\gamma}(B)\}\subset\mathbf{H}^{2(p-1)}_{\gamma}(B)$.


 Эти результаты дополняют  утверждения Теоремы 2b.

\newpage
Реферат:

Ряды Фурье оператора градиент дивергенции 
 и пространства Соболева.

  Р.С.Сакс

Автор изучает структуру пространства
  $\mathbf{L}_{2}(G)$ вектор-функций, квадратично интегрируемых
      по ограниченной односвязной  области $G$ трехмерного    пространства 
     с гладкой границей, 
    и роль операторов   градиента дивергенции     и ротора
     в построении базисов в его ортогональных подпространствах
     ${\mathcal{{A}}}$ и  ${\mathcal{{B}}}$.

В  ${\mathcal{{A}}}$ и  ${\mathcal{{B}}}$ выделяются
подпространства ${\mathcal{{A}}_{\gamma}}(G)\subset{\mathcal{{A}}}$
и $\mathbf{V}^{0}(G)\subset {\mathcal{{B}}}$.
   Оказывается, что  операторы  градиент дивергенции и ротор
продолжаются   в эти подпространства, их расширения $\mathcal{N}_d$
\,и $S$ являются само-сопряженными и обратимыми, а их обратные
операторы $\mathcal{N}_{d}^{-1}$   и $S^{-1}$ -- вполне
непрерывными.

В каждом из этих подпространств  мы строим ортонормированный базис.
Объединяя эти базисы,  получаем полный ортонормированный базис
объемлющего пространства $\mathbf{L}_{2}(G)$, составленный из
собственных функций операторов градиента дивергенции и ротора.

 В случае, когда область $G$ есть шар  $B$,   базисные функции
  определяются элементарными функциями.

Определены пространства $\mathcal{A}^{s}_{\mathcal{K}}(B)$.
Доказано, что условие $\mathbf{v}\in
\mathcal{A}^{s}_{\mathcal{K}}(B)$ необходимо и достаточно для
сходимости ее ряда Фурье (по собственным функциям градиента
дивергенции) в норме пространства Соболева $\mathbf{H}^{s}(B)$.

 Используя ряды Фурье
  функций $\mathbf{f}$ и $\mathbf{u}$,
автор исследует разрешимость (в пространствах $\mathbf{H}^{s}(G)$)
краевой задачи: 
$\nabla\mathbf{div}\mathbf{u}+\lambda\mathbf{u}=\mathbf{f}$ в $G$,
 $\mathbf{n}\cdot\mathbf{u}|_{\Gamma}=g$ на границе,
при условии $\lambda\neq 0$.

В шаре $B$ краевая задача:
$\nabla\mathbf{div}\mathbf{u}+\lambda\mathbf{u}=\mathbf{f}$ в $B$,
  $\mathbf{n}\cdot\mathbf{u}|_S=0$,
 решена полностью и для любых $\lambda$.

Доказано, что при  $\lambda\,\overline{\in}\, Sp\,
(-\nabla\mathrm{div})$
 оператор $\nabla\text{div}+\lambda I$ осуществляет гомеоморфизм
 пространств  $\mathbf{{H}}^2_{\gamma\delta\gamma}(B)$
 и  $\mathbf{{F}}_{\gamma}^0(B)$.

       Перевод

 Fourier series of the $\nabla\,div$  operator and Sobolev  spaces
 II

Saks Romen Semenovich

The author studies structure of space $\mathbf{L}_{2}(G)$of vectors
- functions, which are integrable  with a square of the module on
the bounded domain $G $of three-dimensional space with smooth
boundary, and role of the   gradient of divergence and curl
operators   in construction of bases in its orthogonal subspaces
$\mathcal{A} $and $\mathcal{B} $.

The ${\mathcal{{A}}} $and ${\mathcal{{B}}} $are contain subspaces
${\mathcal{{A}}_{\gamma}} (G) \subset{\mathcal{{A}}} $and
$\mathbf{V}^{0}(G)\subset{\mathcal{{B}}}$.

The gradient of divergence and a curl operators have continuations
in these subspaces, their expansion $\mathcal{N}_d$ and $S$ are
selfadjoint and convertible,and their inverse operators
$\mathcal{N}_{d} ^{-1} $and $S ^{-1} $are compact.

In each of these subspaces we build ortonormal basis. Uniting these
bases, we receive complete ortonormal basis of whole
 space $\mathbf{L}_{2} (G) $, made from eigenfunctions
of the gradient of divergence and curl operators .

In a case, when the domain $G $is a ball $B $, basic functions are
defined by elementary functions.

The spaces $ \mathcal {A} ^ {s} _ {\mathcal {K}} (B) $ are  defined.
Is proved, that condition $ \mathbf {v} \in\mathcal {A} ^ {s} _
{\mathcal {K}} (B) $ is necessary and sufficient for convergence of
its Fourier series  (on eigenfunctions of a gradient of divergence)
in norm of  Sobolev space $ \mathbf {H} ^ {s} (B) $.

Using Fourier series of
  functions $ \mathbf {f} $ and $ \mathbf {u} $,
the author investigates solvability (in spaces $ \mathbf {H} ^ {s}
(G) $) boundary value problem: $ \nabla\mathbf {div}\mathbf {u} +
\lambda\mathbf {u} = \mathbf {f} $ in $G $,  $ \mathbf {n}
\cdot\mathbf {u} | _ {\Gamma} =g $ on  boundary,
 under condition of $ \lambda\neq 0$.

In a ball $B $ a boundary value problem: $ \nabla\mathbf {div}
\mathbf {u} + \lambda\mathbf {u} = \mathbf {f} $ in $B$,
 $ \mathbf {n} \cdot\mathbf {u} | _S=0 $,
 is solved completely and for any $ \lambda $.

It's proved, that at $ \lambda \,\overline {\in} \, Sp \, (
-\nabla\mathrm {div}) $
 the operator $ \nabla\text {div} + \lambda I $ carries out
gomeomorphism (g.e.
  one to one and mutually   continuous mapping)
 of spaces $ \mathbf {{H}} ^2 _ {\gamma\delta\gamma} (B) $ and
 $ \mathbf {{F}} _ {\gamma} ^0 (B) $.

Сакс Ромэн Семенович ведущий научный сотрудник
 Институт Математики с
ВЦ УНЦ РАН 450077, г. Уфа, ул. Чернышевского, д.112 телефон: (347)
272-59-36
                 (347) 273-34-12
факс:        (347) 272-59-36 телефон дом.: (347) 273-84-69 моб.
+79173797538

 e-mail: romen-saks@yandex.ru

 \end{document}